\documentclass[a4paper, BCOR=0.0mm, DIV=calc]{scrartcl}
\usepackage[T1]{fontenc}
\usepackage[utf8]{inputenc}
\usepackage[ngerman]{babel}
\usepackage[osf,sc]{mathpazo}
\usepackage{textcomp}
\usepackage{microtype}
\usepackage{amsmath, amssymb, amsfonts}
\usepackage{tabularx}
\KOMAoptions{twoside=false, twocolumn=false, headinclude=false, footinclude=false, mpinclude=false, pagesize=auto}
\recalctypearea

\newcommand{\etc}{\text{etc}} 
\setkomafont{section}{\mdseries\scshape\Large}
\setkomafont{title}{\mdseries\scshape}

\begin{document}
\title{On divergent series\footnote{
Original title: "`De seriebus divergentibus"', first published in "`\textit{Novi Commentarii academiae scientiarum Petropolitanae} 5, 1760, pp. 205-237"', reprinted in "`\textit{Opera Omnia}: Series 1, Volume 14, pp. 585 - 617"', Eneström-Number E247, translated by: Alexander Aycock}}
\author{Leonhard Euler}
\date{}
\maketitle
\paragraph{§1} 

Because convergent series are defined in that manner, that they consist of continuously decreasing terms, that finally, if the series continues to infinity, vanish completely; it is easily seen, that those series, whose infinitesimal terms do not become nothing, but either stay finite or grow to infinty, have, because they are not convergent, to be referred to the class of divergent series. Depending on whether the last terms of the series, to which one gets in the progression continued to infinity, are either of a finite magnitude or infinite, one has two kinds of divergent series, both of which are further subdivided into two subkinds, depending on whether all terms are affected by the same sign, or the signs $+$ and $-$ alternate with one another. Therefore we will in total have four species of divergent series, from which for the sake of greater clarity I want to add some examples.

\begin{alignat*}{9}
	&\text{I.} \quad \quad \quad && 1 &&+ 1 &&+ 1 &&+ 1 &&+ 1 &&+ ~ 1 &&+ \etc. \\
	&                &&\frac{1}{2} &&+ \frac{2}{3} &&+ \frac{3}{4} &&+ \frac{4}{5} &&+ \frac{5}{6} &&+ ~\frac{6}{7} &&+ \etc.\\[2mm]
	&\text{II.} \quad \quad \quad && 1 &&- 1 &&+ 1 &&- 1 &&+ 1 &&- ~1 &&+ \etc. \\
	& &&\frac{1}{2} &&- \frac{2}{3} &&+ \frac{3}{4} &&- \frac{4}{5} &&+ \frac{5}{6} &&- ~\frac{6}{7} &&+ \etc. \\[2mm]
	\end{alignat*}
	\begin{alignat*}{9}
	&\text{III.} \quad \quad \quad && 1 &&+ 2 &&+ 3 &&+ 4 &&+ 5 &&+ ~6 &&+ \etc.\\
	& &&1 &&+ 2 &&+ 4 &&+ 8 &&+ 16 &&+ 32 &&+ \etc. \\[2mm]
	& \text{IV.} \quad \quad \quad && 1 &&- 2 &&+ 3 &&- 4 &&+ 5 &&- ~6  &&+ \etc. \\
	& &&1 &&- 2 &&+ 4 &&- 8 &&+ 16 &&- 32 &&+ \etc. 
\end{alignat*}
\paragraph{§2}
There is a great disagreement about divergent series of this kind between the mathematicians, while some negate, others do not, that they can be comprehended in one sum. And at first it is certainly clear, that the sums of those series, I referred to the first class, are indeed infinite, because by actually collecting the terms one gets to a sum greater than any number: Hence there is no doubt, that sums of series of this kind can be exhibited by expressions like $\frac{a}{0}$. So the great controversy between Geometers is mainly about the remaining three species; and the argument, which are urged by both sides to defend their postion, are so much convincing, that neither party could be forced to agree with the other.
\paragraph{§3}
From the second species \textsc{Leibniz} at first considered this series
\[
	1 - 1 + 1 - 1 + 1 - 1 + 1 - 1  + \etc.,
\]
the sum of which he stated to be $=\frac{1}{2}$, while basing it on these fairly solid arguments: Hence at first this series arises, if this fraction $\frac{1}{1+a}$ by an iterated division in usual manner is resolved into this series
\[
	1 - a + a^2 - a^3 + a^4 - a^5 + \etc.
\] 
and the the value of the letter $a$ is taken equal to the unity. Then indeed, to further confirm this and to persuade those, who are not used to such calculations, he gave the following argument: If this series is terminated at some point, and the number of terms was even, then its value will be $=0$, but if the number of terms is odd, its value will be $=1$: Therefore, if this series proceds to infinity, and the number of terms can neither be seen to be even nor odd, he concluded, that the sum can neither be $=0$ nor $=1$, but has a certain mean value, being equally different from both, which is $=\frac{1}{2}$.
\paragraph{§4}
To these arguments the adversaries used to object, that at first the fraction $\frac{1}{1+a}$ is only equal to the infinite series
\[
	1 - a + a^2 - a^3 + a^4 - a^5 + a^6 - \etc.,
\]
if $a$ is a fraction smaller than the unity. Hence if the division is abrupted anywhere and the correspondig portion from the remainder is added to the quotient, it will lead to wrong results; hence it is
\[
	\frac{1}{1+a} = 1 - a + a^2 - a^3 + a^4 - \dots \pm a^n \mp \frac{a^{n+1}}{1+a},
\]
and if the number $n$ is put to be  infinite, it is nevertheless not possible, to omit the added fraction $\mp\frac{a^{n+1}}{1+a}$, if it does not really vanish, what is only true in the cases, where $a<1$, and then the series converges. But in the remaining cases one has to have regard for this mantissa $\mp \frac{a^{n+1}}{1+a}$, and although it is affected by the ambiguous sign $\mp$, depending on whether $n$ is either even or odd, it can therefore, if $n$ is infinite, not be neglected, because an infinite number is neither even nor odd, and so one has no reason, what sign is to be preferred. Since it is absurd to believe, that there is no wohle number, not even an infinite one, which is neither even or odd.
\paragraph{§5}
But to this objection the ones, that assign certain sums to divergent series, justifiably answer, that an infinite number is treated as a certain number, and is therefore either even or odd, although it is not determined. When a series is said, to go on to infinity, this contradicts the idea, if a certain term of the series is treated as the last or infinitesimal one: And therefore the objection raised before and concering the mantissa, that has to be added, or subtracted, vanishes by itself. Because in an infinite series one never gets to an end, one therefore never reaches such a place, where it would be necessary to add that  mantissa; and hence this mantissa can not only be neglected but also has to, because it is never left space. And these arguments, that are urged either for or against the divergent series, also concern the fourth species, which usually creates no other doubts than the ones mentioned.
\paragraph{§6}
But those, who argue against the sums of divergent series, have the opinion, that the third species provides them with the best arguments. Hence although the terms of these series increase continuously, one can therefore, by actually collecting the terms, get to a sum greater than any assignable number, that is, by definition, infinity, the defenders of sums in this species are nevertheless forced to admit series of such a kind, whose sums are finite, and even negative or smaller than nothing. Because the fraction $\frac{1}{1-a}$, expanded into a series, yields:
\[
	1 + a + a^2 + a^3 + a^4 + \etc.
\]
the following equations would have to hold:
\begin{alignat*}{9}
&-1 &&= 1 &&+ 2 &&+ 4 &&+ ~8 &&+ 16 &&+ \etc. \\
&-\tfrac{1}{2} &&= 1 &&+ 3 &&+ 9 &&+ 27 &&+ 81 &&+ \etc.
\end{alignat*}
what seems, quite understandably, very suspect to adversaries, because by the addition of only affirmative terms one can never get a negative sum. And hence the more they stress the before mentioned mantissa, that has to be added, because, after having added it, it is perspicuous, that it will be
\[
	-1 = 1 + 2 + 4 + 8 + \dots + 2^n + \frac{2^{n+1}}{1-2},
\]
even though $n$ is an infinite number.
\paragraph{§7}
Therefore the defenders of sums of divergent series, to explain this great paradox, rather subtle, than true, state a difference between negative quantities, while on the one hand smaller than nothing, they argue, on the other hand they are graeter than infinity or more than infinite numbers. On the one hand  they have to accept the value of $-1$, whenever it is imagined, that it arises from the subtraction of the greater number $a+1$ from the smaller $a$, but on the other hand, whenever it is found to be equal to the series $1+2+4+8+16+$etc. and emerges from the division of the number $+1$ by the number $-1$; in that case the number is of course smaller than nothing, but in this  one  greater than infinity. For the sake of further confirmation they give this example of the series of fractions
\[
	\frac{1}{4},\quad \frac{1}{3},\quad \frac{1}{2},\quad \frac{1}{1},\quad \frac{1}{0}, \quad \frac{1}{-1},\quad \frac{1}{-2},\quad \frac{1}{-3}\quad \etc.,
\]
that, because in the first terms it is seen to grow, it is also to be seen to grow continuously, whence they conclude, that it will be $\frac{1}{-1}>\frac{1}{0}$ and $\frac{1}{-2}>\frac{1}{-1}$ and so on; and therefore, if $\frac{1}{-1}$ is expressed by $-1$ and $\frac{1}{0}$ by $\infty$, that $-1>\infty$ and even more $\frac{-1}{2}>\infty$; and in this way they quite ingeniously repel the apparent absurdity.
\paragraph{§8}
Although this distinction seems to be an ingenious idea, it is nevertheless hardly satisfactory for the adversaries and hence seems to violate the certitude of analysis. Hence if the two values of one $-1$, if it is either $=1-2$ or $-\frac{1}{-1}$, are indeed different from each other, that they cannot be confounded, the certitude and the application of the rules, that we follow in calculus, are abolished completely, what would certainly be more absurd than that, this distinction was actually made for; but if it is $1-2=\frac{1}{-1}$, as the precepts of algebra postulate, the task is in no way completed, because the quantity $-1$ itself, that is stated to be equal to the series $1+2+4+8+$etc., is nevertheless smaller and the same difficulty remains. But it nevertheless seems to be true, if we say, that the  same quantities, that are smaller than nothing, can at the same moment be seen as greater than infinity. Hence not only from algebra but also from geometry we know, that there is a jump from positive to negative numbers, the one at zero or nothing, the other at infinity, and therefore the quantities form zero, as by increasing as decreasing, will return to themselves and will finally reach the same term $=0$ again, so that the quantities greater than infinity are also smaller than nothing and the quantities smaller than infinity also correspond to the quanities greater than nothing.
\paragraph{§9}
But the same, who negate that these sums of divergent series, which are usually assigned to them, are correct and justified, do not only not proffer other suggestions, but also state, that they totally believe, that the sum of a divergent series is imaginary. The sum of convergent series as this one
\[
	1 + \frac{1}{2} + \frac{1}{4} + \frac{1}{8} + \frac{1}{16}+\frac{1}{32} + \etc.
\]
can only be admitted to be $=2$, because, the more terms of this series we actually add, the closer we get to two; but for divergent series the matter behaves totally different; hence the more terms we add, the more the sums, that arise, differ from each other and they to not get closer to a certain determined value. From this they conclude, that not even the idea of a sum can be transferred to divergent series and the work, that was consumed by investigating the sums of divergent series, of those is completely useless and contrary to the true principles of analysis.
\paragraph{§10}
But although this difference seems to be real, none of the two parties can be convicted of an error by the other, as often as the use of series of this kind occurs in analyis; it has to be of a great ponderosity, that no party made any mistakes, but the whole dissent lies only in the words and formulations. Hence if in a calculation I get to this series $1-1+1-1+1-1+\text{etc.}$ and substitute $\frac{1}{2}$ for it, certainly no one will ascribe an error to me,  that nevertheless would occur to everybody, if I had put another value in the place of the series; hence there can remain no doubt, that the series $1-1+1-1+1-1+\text{etc.}$ and the fraction $\frac{1}{2}$ are equivalent quantities. So the whole question seems to trace back to the one, whether we correctly call the fraction $\frac{1}{2}$ the sum of the series $1-1+1-1+\text{etc.}$; because those persistently deny this, although they do not dare to deny the equivalence, it is to be feared, that they slip into wrong logic.
\paragraph{§11}
But I believe, that the whole dispute can easily be settled, if we pay close attention to the following. As often as in analysis we get to an either rational or transcendental experession, we usually converted it into an appropriate series, to which the following calculation is more conveniently applied. Hence if infinite series occur in analysis, they arose from the expansion of a certain finite expression, and therefore in a calculation it is always possible, to substitute the formula, from whose expansion the series arose, for the series. Hence as with the greatest gain the rules, to convert finite expressions, but of a less suitable form, into infinte series, were given, vice versa the rules, by which, if any infinte series was given, the finite expression can be found, from which it resulted, have to be considered of the greatest use; and because this expression can always without an error be put in place of the infinite series, it is necessary, that the value of both is the same; hence it is caused, that there is no series, that cannot at the same moment be considered to be equivalent to the finite expression. 
\paragraph{§12}
Hence if we just change the usual notion of a sum in such a way, that we say, that the sum of a certain series is the finite expression, from whose expansion that series itself arises, all difficulties, which were mentioned by both parties, will disappear by itself. Hence at first the expression, from which a convergent series arises, at the same moment exhibits its sum, in the usual sense, and if not, if the series was divergent, the question cannot be seen as absurd any longer, if we find the finite expression, that, expanded according to the analytical rules, produces the series itself. And because it is possible to substitute this expression for its series in a calculation, we will not be able to doubt, that they will even be equal to each other. Having explained this we do not even recede from the usual  notion, if we  call the expression, that it is equal to a certain sum, its sum too, as long as we do not combine the notion with the idea of a sum for divergent series, that, the more terms are actually collected, the closer one has to get to the value of the sum.
\paragraph{§13}
Having said all this in advance, I believe that there will be nobody, who thinks, that I have to be reprehended, because I inquire into the sum of the following series more diligently
\[
	1 - 1 + 2 - 6 + 24 - 120 + 720 - 5040 + 40320 - \mathrm{etc.},
\]
which is the, called this way by \textsc{Wallis}, hypergeometric series, just with alternating signs. This series seems noteworthy all the more, because I have tried several summation methods, that were quite heplful for other tasks of this kind, without success here. At first it is certainly possible to doubt, whether this series has a finite sum or not, because it  diverges even more than any divergent series; but that the sum of the geometric series is finite, was clarified. But because for the geometric series the divergence is not an obstacle, that they are summable, it seems probable, that also this hypergeometric series has a finite sum. So one in numbers, at least approximately, looks for the value of that finite expression, from whose expansion the given series itself arises.
\paragraph{§14}
At first I used the method, based on this foundation: If a series of this kind is given
\[
	s = a - b + c - d + e - f + g - h + \mathrm{etc.}
\]
and, after having neglected the signs of the terms $a$, $b$, $c$, $d$, $e$, $f$ etc., one takes the differences
\[
	b-a,\quad c-b,\quad d-c,\quad e-d\quad \mathrm{etc.}
\]
and further their differences
\[
	c-2b+a,\quad d-2c+b,\quad e-2d+c\quad \mathrm{etc.},
\]
which are called the second differences, and in the same way searches the third, fourth, fifth differences etc., then, if the first terms of these first, second, third, fourth differences etc. are $\alpha$, $\beta$, $\gamma$, $\delta$ etc., I say, the sum of the same given series will be
\[
	s = \frac{a}{2} - \frac{\alpha}{4} + \frac{\beta}{8} - \frac{\gamma}{16} + \frac{\delta}{32} - \mathrm{etc.},
\]
which series, if it is not already convergent, will nevertheless be a lot more convergent than the given one; hence, if the same method is then again applied to this last series, the value of the desired sum expressed by $s$ will be found by means of an even more convergent series.
\paragraph{§15}
This method has the greatest use for summing divergent series of the second and the fourth species, whether one finally reaches constant differences or not, as long as the divergence is not too strong: If it is
\[
	s = 1 - 1 + 1 - 1 + 1 - \text{etc.},
\]
because of
\[
	a = 1,\quad \alpha = 0,\quad \beta = 0\quad \text{etc.}
\]
it will be
\[
	s = \frac{1}{2}.
\]\\

If
\[
\begin{array}{rlcccccccccccccc}
	s &=& 1 &-& 2 &+& 3 &-& 4 &+& 5 &-& 6 &+& \text{etc}, \\
	\text{diff\, I.} & & &1& &1& &1& &1& &1&
\end{array}
\]
it will be
\[
	s = \frac{1}{2} - \frac{1}{4} = \frac{1}{4},
\]
as it is known from elsewhere.\\

If it is
\[
\begin{array}{llcccccccccccccc}
	s &=& 1 &-& 4 &+& 9 &-& 16 &+& 25 &-& 36 &+& \text{etc.}, \\
	\text{diff\,  I.} & & &3& &5& &7& &9& &11& \\
	\text{diff\, II.} & & & &2& &2& &2& &2& &
\end{array}	
\]
it will be
\[
	s = \frac{1}{2} - \frac{3}{4} + \frac{2}{8} = 0,
\]
as it is also known. \\

If it is
\[
\begin{array}{llcccccccccccccc}
	s &=& 1 &-& 3 &+& 9 &-& 27 &+& 81 &-& 243 &+& \text{etc.}, \\
	\text{diff\, I.}   & & &2& &6& &18& &54& &162& \\
	\text{diff\, II.}  & & & &4& &12& &36& &108& & \\
	\text{diff\, III.} & & & & &8& &24& &72& & & \\
	\text{diff\, IV.}  & & & & & &16& &48& & & & \\
	                  & & & & & &  & \text{etc.}&  & & & & 
\end{array}	
\]
it will be
\[
	s = \frac{1}{2} - \frac{2}{4} + \frac{4}{8} - \frac{8}{16} + \text{etc.} = \frac{1}{2} - \frac{1}{2} + \frac{1}{2} - \frac{1}{2} + \text{etc.} = \frac{1}{4}
\]
sein.
\paragraph{§16}
Now let us apply this method to the proposed series
\[
	A = 1 - 1 + 2 - 6 + 24 - 120 + 720 - 5040 + 40320 - \etc.,
\]
which because of $1-1 = 0$, if it is divided by $2$, changes into
\begin{center}
$\frac{A}{2} = 1 - 3 + 12 - 60 + 360 - 2520 + 20160 - 181440 + \etc.$\\
$2,\quad 9,\quad 48,\quad 300,\quad 2160,\quad 17640,\quad 161280$\\
$7,\quad 39,\quad 252,\quad 1860,\quad 15480,\quad 143640$ \\
$32,\quad 213,\quad 1608,\quad 13620,\quad 128160$\\
$181,\quad 1395,\quad 12012,\quad 114540$\\
$1214,\quad 10617,\quad 102528$\\
$9403,\quad 91911$\\
$82508$\\
\end{center}
Hence it follows, that it will be
\[
	\frac{A}{2} = \frac{1}{2} - \frac{2}{4} + \frac{7}{8} - \frac{32}{16} + \frac{181}{32} - \frac{1214}{64} + \frac{9403}{128} - \frac{82508}{256} + \etc.
\]
or
\begin{center}
	$A = \dfrac{7}{4} - \dfrac{32}{8} + \dfrac{181}{16} - \dfrac{1214}{32} + \dfrac{9403}{64} - \dfrac{82508}{128} + \etc.$\\[1ex]
	$\dfrac{18}{8},\quad \dfrac{117}{16},\quad \dfrac{852}{32},\quad \dfrac{6975}{64},\quad \dfrac{63702}{128}$\\[1ex]
	$\dfrac{81}{16},\quad \dfrac{618}{32},\quad \dfrac{5271}{64},\quad \dfrac{49752}{128}$\\[1ex]
	$\dfrac{456}{32},\quad \dfrac{4035}{64},\quad \dfrac{39210}{128}$\\[1ex]
	$\dfrac{3123}{64},\quad \dfrac{31140}{128}$\\[1ex]
	$\dfrac{24894}{128}$
\end{center}
Therefore
\[
	A = \frac{7}{8} - \frac{18}{32} + \frac{81}{128} - \frac{456}{512} + \frac{3123}{2048} - \frac{24894}{8192} + \etc.
\]
or
\begin{center}
	$A - \dfrac{5}{16} = \dfrac{81}{128} - \dfrac{456}{512} + \dfrac{3123}{2048} - \dfrac{24894}{8192} + \etc$ \\[1ex]
	$\quad\dfrac{132}{512},\quad \dfrac{1299}{2048},\quad \dfrac{12402}{8192}$\\[1ex]
	$\quad\dfrac{771}{2048},\quad \dfrac{7206}{8192}$\\[1ex]
	$\quad\dfrac{4122}{8192}$
\end{center}
So
\[
	A - \frac{5}{16} = \frac{81}{256} - \frac{132}{2048} + \frac{771}{16384} - \frac{4122}{131072}
\]
or
\[
	A = \frac{5}{16} + \frac{516}{2048} + \frac{2046}{131072} + \etc = \frac{38015}{65536} = 0,580.
\]
Hence it is clear, that the sum of this series is nearly $=0,580$; but because of the neglected terms it will be a little bit greater, what agrees very well with the things, that are to be demonstrated below, where the sum of this series will be shown to be $=0,59634736$; at the same moment it is indeed clear, that this method is apt enough, to find the sum that exact.
\paragraph{§17}
Next I tried it this way: Let this series be given
\[
\begin{array}{llllllllllllllll}
	& 1 & 2 & 3 & ~4 & ~5 & ~~6 & ~~~7 & \dots & n & ~n+1 \\
	B)  & 1, & 2, & 5, & 16, & 65, & 326, & 1957, & \dots & P, & nP+1 
\end{array}
\]
the differences are
\begin{center}
$1,\quad 3,\quad 11,\quad 49,\quad 261,\quad 1631$ \\
$2,\quad 8,\quad 38,\quad 212,\quad 1370$ \\
$6,\quad 30,\quad 174,\quad 1158$ \\
$24,\quad 144,\quad 984$ \\
$120,\quad 840$ \\
$720$
\end{center}
because the first terms of its continued differences are
\[
	1,\quad 2,\quad 6,\quad 24,\quad 120,\quad 720\quad \etc.,
\]
the term corresponding to the exponent $n$ will be
\begin{align*}
	P =& 1 + (n-1) + (n-1)(n-2) + (n-1)(n-2)(n-3) \\
	&+ (n-1)(n-2)(n-3)(n-4) + \etc.
\end{align*}
Since, if $n=0$, the term corresponding to the exponent $0$ or preceding the first will be
\[
	1 - 1 + 2 - 6 + 24 - 120 - \etc. = A,
\]
so that, if the term corresponding to the exponent $0$ of this series could be found, the same simultaneously would be the value or the sum of the given series
\[
	A = 1 - 1 + 2 - 6 + 24 - 120 + 720 - \etc.
\]
Hence if that series $B$ is inverted, that one has the series
\begin{alignat*}{19}
	&  1 \quad &&  2 \quad &&  3 \quad && ~4 \quad && ~5 \quad && ~~6 \quad && ~~~7 \\
	C) \quad  & 1, \quad && \frac{1}{2}, \quad && \frac{1}{5}, \quad && \frac{1}{16}, \quad  && \frac{1}{65}, \quad && \frac{1}{326}, \quad && \frac{1}{1957} \quad  && \etc .
\end{alignat*}
the term corresponding to the exponent $0$ of this series will be $=\frac{1}{A}$, whence the value of $A$ can be perceived from it. Let the single differences of this series begin with the terms $\alpha$, $\beta$, $\gamma$, $\delta$, $\varepsilon$ etc., of course by taking the difference in such a way, that any term is subtracted from the preceding; the term corresponding to the exponent $n$ will be
\[
	\frac{1}{P} = 1 - (n-1)\alpha + \frac{(n-1)(n-2)}{1\cdot 2}\beta - \frac{(n-1)(n-2)(n-3)}{1\cdot 2\cdot 3}\gamma + \etc.
\]
und hence for $n=0$ it will be by means of a surely converging series
\[
	\frac{1}{A} = 1 + \alpha + \beta + \gamma + \delta + \etc.
\]
It is indeed, by converting these fractions into decimals,
\[
\begin{array}{rccccc}
& \text{diff. } 1 & \text{diff. } 2 & \text{diff.} 3 & \text{diff. } 4 & \text{diff. }  5 \\
1 = 1,0000000 & & & & & \\
& 5000000 & & & & \\
\frac{1}{2} = 0,5000000 &  & 2000000 & & & \\ 
& 3000000 & & 375000 & & \\
\frac{1}{5} = 0,2000000 & & 1625000 & & -346154 & \\
& 1375000 & & 721154 & & -511445 \\
\frac{1}{16} = 0,0625000 & & 903848 & & +165291 & \\
& 471154 & & 555863 & & -140195\\
\frac{1}{65} = 0,0153846 & & 347983 & & +305486 & \\
& 123171 & & 250377 & & +131530\\
\frac{1}{326} = 0,0030675 & & 97606 & & +173956 & \\
& 25565 & & 76421 & & +114979\\
\frac{1}{1957} =0,0005110 & & 21185 & & +58977 & \\
& 4380 & & 17444 & & +44716\\
0,0000370 & & 3741 & & +14261 & \\
& 639 & & 3183 & & +11564\\
0,0000091 & & 558 & & +2697 & \\
& 81 & & 486 & & +2275\\
0,0000010 & & 72 & & +422 & \\
& 9 & & 64 & & +365\\
0,0000001 & & 8 & & +57 &  
\end{array}
\]
From this differences it will therefore be
\[
	\frac{1}{A} = 1,6517401 \quad \text{and} \quad A = 0,6,
\]
which argees to with the value found before to a high enough degree of accuracy; but because of the fourth, fifth and some of the following differences this method is nevertheless not certain enough.
\paragraph{§18}
Let us take the logarithms of the single terms of the series $B$, that one has this new series
\[
\begin{array}{ccccccccccl}
	& 1 & 2 & 3 & 4 & 5 & 6 & 7 & 8 \\
	D) & \log{1}, &\log{2}, & \log{5}, & \log{16}, & \log{65}, & \log{326}, & \log{1957}, & \log{13700}& \etc. 
\end{array}
\]
in whose continued differences taken the usual way the first terms shall be $\alpha$, $\beta$, $\gamma$, $\delta$, $\varepsilon$ etc., and the term corresponding to the exponent $0$ of this series will be
\[
	0 - \alpha + \beta - \gamma + \delta - \varepsilon + \etc.,
\]
which will therefore will be the logarithm of the desired sum $=A$. The logarithms with the continued differences are indeed the following:
\[
\begin{array}{ccccccccc}
& \text{diff. } 1 & \text{diff. }2 & \text{diff. }3 & \text{diff. }4 & \text{diff. }5 & \text{diff. }6 & \text{diff. } 7 & \text{diff. } 8 \\
0,0000000 & & & & & & & & \\
& 0,3010300 & & & & & & & \\
0,3010300 & & 969100 & & & & & & \\
& 0,3979400 & & 103000 & & & & & \\
0,6989700 & & 1072100 & & -138666& & & & \\
& 0,5051500 & & -35666 & & +53006& & & \\
1,2041200 & & 1036434 & & -85660 & & +19562 & & \\
& 0,6087934 & & -121326 & & +72568 & & -57744 & \\
1,8129134 & & 915108 & & -12092 & & -38182 & & +65446 \\
& 0,7003042 & & -134418 & & +34386 & & +7702 & \\
2,5132176 & & 780690 & & +21294 & & -30480 & & \\
& 0,7783732 & & -113124 & & +3906 & & & \\
3,2915908 & & 667566 & & +25200 & & & & \\
& 0,8451298 & & -87925 & & & & & \\
4,1367206 & & 579641 & & & & & & \\
& 0,9030939 & & & & & & & \\
5,0398145 & & & & & & & &
\end{array}
\]
hence it will be
\[
\begin{array}{rcccccc}
& \text{diff. }1 & \text{diff. } 2 & \text{diff. }3 & \text{diff. }4 & \text{diff. }5 & \text{diff. }6 \\  
\log{A} =-0,3010300 & & & & & &\\
& +2041200 & & & & & \\
+0,0969100 & & +1175100 & & & & \\
& +866100 & & +550666 & & & \\
-0,0103000 & & +624434 & & +359570 & & \\
& +241666 & & +191096 & & +826928 & \\
-0,0138666 & & +433338 & & -467358 & & +2133994 \\
& -191672 & & +658454 & & -1307066 & \\
-0,0053006 & & -225116 & & +839708 & & -2083670\\
& +33444 & & -181254 & & +776604 & \\
+0,0019562 & & -43862 & & +63103 & & \\
& +77306 & & -244357 & & & \\
+0,0057744 & & +200495 & & & & \\
& -123189 & & & & & \\
+0,0065445 & & & & & &
\end{array}
\]
whence by the method explained it will be 
\[
	\log{\frac{1}{A}} = \frac{0,0310300}{2} + \frac{2041200}{4} + \frac{1175100}{8} + \frac{550666}{16} + \frac{359570}{32} + \frac{826928}{64} + \etc.
\]
or
\[
	\log{\frac{A}{1}} = 0,7779089 \quad \text{and therefore} \quad A = 0,59966,
\]
which number may easily be calculated to be still greater than the true one. Nevertheless even on this way one can neither certain enough nor comfortable enough get cognition of the value $A$, even though this method yields an infinite amount of ways to investigate this value; but from those the ones certainly seem much more apt for this purpose than others.
\paragraph{§19}
Now let us also investigate the value of this series analytically, but let us accept it in a broader sense; it shall be
\[
	s = x - 1x^2 + 2x^3 - 6x^4 + 24x^5 - 120x^5 + \etc.,
\]
which differentiated will give
\[
	\frac{\mathrm{d}s}{\mathrm{d}x} = 1 - 2x + 6xx - 24x^3 + 120x^4 - \etc.  = \frac{x-s}{xx},
\]
whence it becomes
\[
	\mathrm{d}s + \frac{s\mathrm{d}x}{xx} = \frac{\mathrm{d}x}{x},
\]
the integral of which equation, if $e$ is taken for that number, whose hyperbolic logarithm is $=1$, will be
\[
	e^{-1:x}s = \int{\frac{e^{-1:x}}{x}\mathrm{d}x} \quad \text{and} \quad s = e^{1:x}\int{\frac{e^{-1:x}}{x}\mathrm{d}x}.
\]
In the case $x=1$ it will be
\[
	1 - 1 + 2 - 6 + 24 - 120 + \etc. = e\int{\frac{e^{-1:x}}{x}\mathrm{d}x}.
\]
Hence this series expresses the area of the curved line, whose nature between the abscissa $x$ and $y$ is contained in this equation
\[
	y = \frac{e\cdot e^{-1:x}}{x},
\]
if the abscissa $x$ is put $=1$, or it will be
\[
	y = \frac{e}{e^{1:x}\cdot x}.
\]
But this curve is conditioned in such a way, that for $x=0$ $y$ becomes $=0$; but if $x=1$, $y=1$; but the intermediate values of the ordinate will indeed behave like this, that
\[
\begin{array}{lcllcl}
\text{if it was} & \qquad \qquad & \text{it will then also be} \qquad & \text{if it was} & \qquad \qquad & \text{it will then also be} \\[3mm]
x = \dfrac{0}{10} & \qquad \qquad & y = 0                 \qquad & x = \dfrac{5}{10} & \qquad \qquad & y = \dfrac{10}{5e^{5:5}} = \dfrac{2}{e} \\[3mm]
x = \dfrac{1}{10} & \qquad \qquad & y = \dfrac{10}{e^{9:1}} \qquad & x = \dfrac{6}{10} & \qquad \qquad & y = \dfrac{10}{6e^{4:6}} \\[3mm]
x = \dfrac{2}{10} & \qquad \qquad & y = \dfrac{10}{2e^{8:2}} \qquad & x = \dfrac{7}{10} & \qquad \qquad & y = \dfrac{10}{7e^{3:7}} \\[3mm]
x = \dfrac{3}{10} & \qquad \qquad & y = \dfrac{10}{3e^{7:3}} \qquad & x = \dfrac{8}{10} & \qquad \qquad & y = \dfrac{10}{8e^{2:8}} \\[3mm]
x = \dfrac{4}{10} & \qquad \qquad & y = \dfrac{10}{4e^{6:4}} \qquad & x = \dfrac{9}{10} & \qquad \qquad & y = \dfrac{10}{9e^{1:9}}
\end{array}
\]
Hence having constructed this curve, it will instantaneously become clear, that its area corresponding to the abscissa $x=1$ is not only finite, but also smaller than the area of the unit square, namely $=1$, but greater than its half $=\frac{1}{2}$. Hence if the base $x=1$ is divided up into ten equal parts and the portions of the area are considered as trapeziods and those areas are investigated, one will obtain this value very close to the true one of the series 
\[
	1 - 1 + 2 - 6 + 24 - 120 + \etc. = A
\] 
namely
\[
	A = 0 + \frac{1}{e^{9:1}} + \frac{1}{2e^{8:2}} + \frac{1}{3e^{7:3}} + \frac{1}{4e^{6:4}} + \frac{1}{5e^{5:5}} + \frac{1}{6e^{4:6}} + \frac{1}{7e^{3:7}} + \frac{1}{8e^{2:8}} + \frac{1}{9e^{1:9}} + \frac{1}{20}.
\]
These terms, because $e=2,718281828$, attain the following values:
\begin{align*}
	\frac{1}{e^{9:1}} &= 0,00012341 \\
	\frac{1}{2e^{8:2}} &= 0,00915782 \\
	\frac{1}{3e^{7:3}} &= 0,03232399 \\
	\frac{1}{4e^{6:4}} &= 0,05578254 \\
	\frac{1}{5e^{5:5}} &= 0,07357589 \\
	\frac{1}{6e^{4:6}} &= 0,08556952 \\
	\frac{1}{7e^{3:7}} &= 0,09306272 \\
	\frac{1}{8e^{2:8}} &= 0,09735007 \\
	\frac{1}{9e^{1:9}} &= 0,09942659 \\
	\frac{1}{20} &= 0,05000000 \\
	\text{hence}\quad A &= 0,59637255 
\end{align*}
which value differs from the true one already in a hardly noticeable way. But if the abscissa would have divided up into more parts, then this value would have been found more precisely.
\paragraph{§20}
Because the sum was found as
\[
	A = \int{\frac{e^{1 - 1:x}}{x}\mathrm{d}x},
\]
let us set
\[
	v = e^{1 - 1:x},
\]
so that for $x=0$ it also is $v=0$ and for $x=1$ it is  $v=1$; it will be $1 - \frac{1}{x} = \log{v}$ and $x = \frac{1}{1-\log{v}}$ and $\log{x} = -\log{(1-\log{v})}$, whence it becomes
\[
	\frac{\mathrm{d}x}{x} = \frac{\mathrm{d}v}{v(1-\log{v})}
\]
Because it is
\[
	A = \int{\frac{v\mathrm{d}x}{x}},
\]
after having set $x=1$ and $v=1$ it will also be
\[
	A = \int{\frac{\mathrm{d}v}{1 - \log{v}}},
\]
having put $v=1$ after the integration. But it will be by integrating by a series term by term
\begin{align*}
	A &= \int{\frac{\mathrm{d}v}{1-\log{v}}} = \frac{v}{1-\log{v}} - \frac{1\cdot v}{(1-\log{v})^2} + \frac{1\cdot 2\cdot v}{(1-\log{v})^3} \\
	&- \frac{1\cdot 2\cdot 3\cdot v}{(1-\log{v})^4} + \frac{1\cdot 2\cdot 3\cdot 4\cdot v}{(1-\log{v})^5} - \etc.
\end{align*}
and for $v=1$ because of $\log{v} = 0$, as we assumed, it will be
\[
	A = 1 - 1 + 1\cdot 2 - 1\cdot 2\cdot 3 + 1\cdot 2\cdot 3\cdot 4 - 1\cdot 2\cdot 3\cdot 4\cdot 5 + \etc.
\]
Hence $A$ will again be the area of the curve, whose nature between the abscissa $v$ and the ordinate $y$ is expressed by this equation
\[
	y = \frac{1}{1- \log{v}},
\]
if the abscissa $v$ is set $=1$, of course, in which case also $y=1$. But it has to be noted, that $\log v$ denotes the hyperbolic logarithm of $v$. Hence having divided the abscissa $v=1$ up into ten parts again, and the ordinates in the single points of the division will behave in this way:
\begin{alignat*}{9}
&\text{if $v$ is} && \qquad &&  \quad \quad \text{$y$ will be} && \text{if $v$ is} && \qquad && \quad \quad \text{$y$ will be} \\[2mm]
&v = \frac{0}{10}, && \qquad \qquad && y = 0; \qquad \qquad \quad &&v = \frac{5}{10}, && \qquad \qquad && y = \frac{1}{(1 + \log{10} - \log{5})};  \\
&v = \frac{1}{10}, && \qquad \qquad && y = \frac{1}{(1 + \log{10} - \log{1})}; \qquad \qquad \quad &&v = \frac{6}{10}, && \qquad \qquad && y = \frac{1}{(1 + \log{10} - \log{6})};\\
&v = \frac{2}{10}, && \qquad \qquad && y = \frac{1}{(1 + \log{10} - \log{2})}; \qquad \qquad \quad &&v = \frac{7}{10}, && \qquad \qquad && y = \frac{1}{(1 + \log{10} - \log{7})};\\
&v = \frac{3}{10}, && \qquad \qquad && y = \frac{1}{(1 + \log{10} - \log{3})}; \qquad \qquad \quad &&v = \frac{8}{10}, && \qquad \qquad && y = \frac{1}{(1 + \log{10} - \log{8})};\\
&v = \frac{4}{10}, && \qquad \qquad && y = \frac{1}{(1 + \log{10} - \log{4})}; \qquad \qquad \quad &&v = \frac{9}{10}, && \qquad \qquad && y = \frac{1}{(1 + \log{10} - \log{9})};\\
&v = \frac{5}{10}, && \qquad \qquad && y = \frac{1}{(1 + \log{10} - \log{5})}; \qquad \qquad \quad &&v = \frac{10}{10}, && \qquad \qquad && y = 1.\\
\end{alignat*}
And therefore by approximation of the area one will again obtain the value of the letter $A$ to a high enough degree of accuracy.
\paragraph{§21}
But there is another method, derived from the nature of continued fractions, to inquire into the sum of this series, which completes the task a lot easier and faster; hence let, by the expressing the formula more generally, be
\[
	A = 1 - 1x + 2x^2 - 6x^3 + 24x^4 - 120x^5 + 720x^6 - 5040x^7 + \etc. = \frac{1}{1+B};
\] 
it will be
\[
	B = \frac{1x - 2x^2 + 6x^3 - 24x^4 + 120x^5 - 720x^6 + 5040x^7 - \etc.}{1 - 1x + 2x^2 - 6x^3 + 24x^4 - 120x^5 + 720x^6 - 5040x^7 + \etc.} = \frac{x}{1+C}
\]
and
\[
	1 + C = \frac{1 - 1x + 2x^2 - 6x^3 + 24x^4 - 120x^5 + 720x^6 - 5040x^7 + \etc.}{1 - 2x + 6x^2 - 24x^3 + 120x^4 - 720x^5 + 5040x^6 - \etc.}.
\]
Therefore
\[
	C = \frac{x - 4x^2 + 18x^3 - 96x^4 + 600x^5 - 4320x^6 + \etc.}{1 - 2x + 6x^2 - 24x^3 + 120x^4 - 720x^5 + \etc.} = \frac{x}{1+D}
\]
hence
\[
	D = \frac{2x - 12x^2 + 72x^3 - 480x^4 + 3600x^5 - \etc.}{1 - 4x + 18x^2 - 96x^3 + 600x^4 - \etc.} = \frac{2x}{1+E}
\]
Further
\[
	E = \frac{2x - 18x^2 + 144x^3 - 1200x^4 + \etc.}{1 - 6x + 36x^2 - 240x^3 + \etc.} = \frac{2x}{1-F}
\]
and
\[
	F = \frac{3x - 36x^2 + 360x^3 - \etc.}{1 - 9x + 72x^2 - 600x^3 + \etc.} = \frac{3x}{1+G}.
\]
It will be
\[
	G = \frac{3x - 48x^2 + \etc.}{1 - 12x + 120x^2 - \etc.}  =\frac{3x}{1+H}.
\]
So
\[
	H = \frac{4x - \etc}{1 - 16x + \etc} = \frac{4x}{1 + I}.
\]
And therefore it will become clear, that it will analogously be
\[
	I = \frac{4x}{1 + K},\quad K = \frac{5x}{1+L},\quad L = \frac{5x}{1+M}\quad \text{etc. to infinity},
\]
so that the structure of these formulas is easily perceived. Having substituted these values one after another it will be
\[
	1 - 1x + 2x^2 - 6x^3 + 24x^4 - 120x^5 + 720x^6 - 5040x^7 + \etc.
\]
\[
	A = \cfrac{1}{1 + \cfrac{x}{1 + \cfrac{x}{1 + \cfrac{2x}{1 + \cfrac{2x}{1 + \cfrac{3x}{1 + \cfrac{3x}{1 + \cfrac{4x}{1 + \cfrac{4x}{1 + \cfrac{5x}{1 + \cfrac{5x}{1 + \cfrac{6x}{1 + \cfrac{6x}{1 + \cfrac{7x}{\etc.}}}}}}}}}}}}}}
\]
\paragraph{§22}
But how the value of continued fractions of this kind are to be investigated, I showed elsewhere. Because the integer parts of the single denominators are unities of course, only the numerators are important for the calculation; hence let $x=1$ and the investigation of the sum $A$ will be performed as follows:
\begin{alignat*}{15}
& \qquad \qquad A ~= ~&&\frac{0}{1}, \quad &&\frac{1}{1}, \quad  &&\frac{1}{2}, \quad &&\frac{2}{3}, \quad &&\frac{4}{7}, \quad &&\frac{8}{13}, \quad &&\frac{20}{34}, \quad &&\frac{44}{73}, \quad &&\frac{124}{209}, \quad &&\frac{300}{501} \quad &&\etc. \\[1ex]
&\text{Numerators}: \quad  &&1, &&1, &&2, &&2, &&3, &&~3, && ~4, &&~4, &&~~5, &&~~5 &&\etc.
\end{alignat*}
The fractions, exhibited here, get continuously closer to the true value of $A$ of course and they are alternately too great and too small, so that it is
\begin{alignat*}{9}
	&A > \frac{0}{1},\quad && A > \frac{1}{2},\quad && A > \frac{4}{7},\quad && A > \frac{20}{34},\quad && A > \frac{124}{209}\quad && \etc. \\
	&A < \frac{1}{1},\quad && A < \frac{2}{3},\quad && A < \frac{8}{13},\quad && A < \frac{44}{73},\quad && A < \frac{300}{501}\quad && \etc.
\end{alignat*}
Hence the values of $A$ will be in decimal numbers
\[
\begin{array}{ccc}
\text{too small values} & \qquad & \text{too great values} \\[2mm]
0,0000000000 & \qquad \qquad & 1,0000000000 \\
0,5000000000 & \qquad \qquad & 0,6666666667 \\
0,5714285714 & \qquad \qquad & 0,6153846154 \\
0,5882352941 & \qquad \qquad & 0,6027397260 \\
0,5933001436 & \qquad \qquad & 0,5988023952 
\end{array}
\]
If now between the too great and too small values, that are respectively next to each other, the arithemtical mean is taken, there will anew emerge alternately too great and too small values, which are the following:
\[
\begin{array}{ccc}
\text{too small values} & \qquad & \text{to great values} \\[2mm]
0,5000000000 & \qquad \qquad & 0,7500000000 \\
0,5833333333 & \qquad \qquad & 0,6190476190 \\
0,5934065934 & \qquad \qquad & 0,6018099548 \\
0,5954875100 & \qquad \qquad & 0,5980205807 \\
0,5960519153 & \qquad \qquad &  
\end{array}
\]
and so we already get quite close to the true value of $A$.
\paragraph{§23}
But we will be able to investigate the value of this fraction part by part in this way: Let 
\[
	A = \cfrac{1}{1 + \cfrac{1}{1+ \cfrac{1}{1 + \cfrac{2}{1 + \cfrac{2}{1 + \cfrac{3}{1 + \cfrac{3}{1 + \cfrac{4}{1 + \cfrac{4}{1 + \cfrac{5}{1 + \cfrac{5}{1 + \cfrac{6}{1 + \cfrac{6}{1 + \cfrac{7}{1 + \cfrac{7}{1 + \cfrac{8}{1 + \cfrac{8}{1 + p}}}}}}}}}}}}}}}}}
\]
and
\[
	p = \cfrac{9}{1 + \cfrac{9}{1 + \cfrac{10}{1 + \cfrac{10}{1 + \cfrac{11}{1 + \cfrac{11}{1 + \cfrac{12}{1 + \cfrac{12}{1 + \cfrac{13}{1+ \cfrac{13}{1 + \cfrac{14}{1 + \cfrac{14}{1 + \cfrac{15}{1 + \cfrac{15}{1 +q}}}}}}}}}}}}}}
\]
and
\[
	q = \cfrac{16}{1 + \cfrac{16}{1 + \cfrac{17}{1 + \cfrac{17}{1 + \cfrac{18}{1+ \cfrac{18}{1 + \cfrac{19}{1 + \cfrac{19}{1 + \cfrac{20}{1 + \cfrac{20}{1 + r}}}}}}}}}}
\]
it will be
\[
	r = \cfrac{21}{1 + \cfrac{21}{1 + \cfrac{22}{1 + \cfrac{22}{1 + \cfrac{23}{1 + \cfrac{23}{1 + \etc.}}}}}}
\]
Having expanded these values, one will at first find
\[
	A = \frac{491459820 + 139931620p}{824073141 + 234662231p},
\]
then
\[
	p = \frac{2381951 + 649286q}{887640 + 187440q}
\]
and
\[
	q = \frac{11437136 + 2924816r}{3697925 + 643025r}.
\]
Hence it remains, that the value of $r$ is defined, what is certainly as difficult as the one of $A$, but it suffices, to know the value of $r$ only approximately here; since a certain error, committed in  the value of $r$, results in a much smaller error in the value of $q$ and hence again causes a lot smaller error in the value of $p$; from this the error, staining the value of $A$, will be completely imperceptible in the end.
\paragraph{§24}
Because further the numerators $21$, $21$, $22$, $22$, $23$ etc. that are included in the continued fraction of $r$, already get closer to the ratio of equality, at least from the beginning, one can obtain help from this to recognize its value. Hence if all these numerators were equal, that it was
\[
	r = \cfrac{21}{1 + \cfrac{21}{1 + \cfrac{21}{1 + \cfrac{21}{1 + \etc.}}}},
\]
it would be
\[
	r = \frac{21}{1+r}
\]
and hence
\[
	rr+r = 21
\]
and 
\[
	r = \frac{\sqrt{85}-1}{2}.
\]
But because these denominators grow, this value will in fact be smaller. Nevertheless it is possible to conclude, if three continued fractions following each other are set to be
\begin{minipage}{0.5\textwidth}
\centering
\[
	r = \cfrac{21}{1 + \cfrac{21}{1 + \cfrac{22}{1 + \cfrac{22}{1 + \cfrac{23}{1 + \etc.}}}}}
\]
\end{minipage}
\begin{minipage}{0.5\textwidth}
\centering
\[
	s = \cfrac{22}{1 + \cfrac{22}{1 + \cfrac{23}{1 + \cfrac{23}{1 + \cfrac{24}{1 + \etc.}}}}}
\]
\end{minipage}
\begin{center}
\[
	t = \cfrac{23}{1 + \cfrac{23}{1 + \cfrac{24}{1 + \cfrac{24}{1 + \cfrac{25}{1 + \etc.}}}}}
\]
\end{center}
that the values of the quantities $r$, $s$, $t$ will proced in an arithmetical progression and it will be $r+t=2s$; hence the value of $r$ will be calculated to a high enough degree of accuracy. But to extend this investigation even further, let us take for the number $21$, $22$, $23$ this indefinite ones $a-1$, $a$ and $a+1$, that it is

\begin{minipage}{0.5\textwidth}
\centering
\[
	r = \cfrac{a-1}{1 + \cfrac{a-1}{1 + \cfrac{a}{1 + \cfrac{a}{1 + \cfrac{a+1}{1 + \etc.}}}}}
\]
\end{minipage}
\begin{minipage}{0.5\textwidth}
\centering
\[
	s = \cfrac{a}{1 + \cfrac{a}{1 + \cfrac{a+1}{1 + \cfrac{a+1}{1 + \cfrac{a+2}{1 + \etc.}}}}}
\]
\end{minipage}
\begin{center}
\[
	t = \cfrac{a+1}{1 + \cfrac{a+1}{1 + \cfrac{a+2}{1 + \cfrac{a+2}{1 + \cfrac{a+3}{1 + \etc.}}}}}
\]
\end{center}
and it will be
\[
	r = \cfrac{a-1}{1 + \cfrac{a-1}{1+s}} \qquad s = \cfrac{a}{1 + \cfrac{a}{1+t}},
\]
whence it is effected
\[
	r = \frac{(a-1)s + a -1 }{s+a}
\]
and
\[
	s = \frac{at+a}{t+a+1} \quad \text{or} \quad t = \frac{(a+1)s-a}{a-s},
\]
whence it becomes
\[
	r+t = \frac{2ss + (2aa - 2a + 1)s - a}{aa - ss} = 2s;
\]
and therefore it will be
\[
	2s^3 + 2ss - (2a-1)s - a = 0,
\]
from which equation one may determine the value of $s$ and further the value of $r$.
\paragraph{§25}
Now let $a=22$ and we will have to solve this cubic equation
\[
	2s^3 + 2ss - 43s - 22 = 0,
\]
whose root is immemidiately discovered to lie beweteen the limits $4$ and $5$. Hence let $s$ be $=4+u$ and it will be
\[
	34 = 69u + 26uu + 2u^3.
\]
Let further $u$ be $=0,4+v$. It will be
\[
	u^2 = 0,16 + 0,8v + vv \quad \text{and} \quad u^3 = 0,064 + 0,48v + 1,2v^2 + v^3
\]
and hence
\[
	2,112 = 90,76v + 28,4v^2 + 2v^3,
\]
hence it will be approximately
\[
	v = 0,023 \quad \text{and} \quad s = 4,423.
\]
Because it is
\[
	r = \frac{21s+21}{s+2},
\]
it will be
\[
	r = \frac{113,883}{26,423} = 4,31
\]
and hence further
\[
	q = \frac{24043093}{6469363} = 3,71645446,
\]
whence one obtains
\[
	p = \frac{4794992,85}{1584252,22} = 3,0266600163
\]
and from this finally
\[
	A = \frac{914985259,27}{1534315932,90} = 0,5963473621372,
\]
which value, converted into a continued fraction, yields
\[
	A = \cfrac{1}{1 + \cfrac{1}{1+ \cfrac{1}{2 + \cfrac{1}{10 + \cfrac{1}{1 + \cfrac{1}{1 + \cfrac{1}{4 + \cfrac{1}{2 + \cfrac{1}{2 + \cfrac{1}{13 + \cfrac{1}{4 + \etc.}}}}}}}}}}},
\]
whence the following values, exhibiting the value of $A$ approximately, are found
\begin{alignat*}{15}
	& && 1 && 1 && 2 && 10 && ~1 && ~1 && ~~4 && ~~2 && ~~~2 && ~~13\\
	&A =~~ && \frac{0}{1}, \quad && \frac{1}{1}, \quad && \frac{1}{2}, \quad && ~\frac{3}{5}, \quad && \frac{31}{52}, \quad && \frac{34}{57}, \quad && \frac{65}{109}, \quad && \frac{294}{493}, \quad && \frac{653}{1095}, \quad && \frac{1600}{2683} \quad && \etc.
\end{alignat*}
But these fraction are alternately greater and smaller than the value of $A$ and the last $\frac{1600}{2683}$ is certainly too large, the excess is nevertheless smaller than $\frac{1}{2683 \cdot 35974}$; hence, because it is
\[
	\frac{1}{A} = \frac{2683}{1600},
\]
it will approximately be
\[
	\frac{1}{A} = 1,676875.
\]
\paragraph{§26}
The method, I used above in §$21$ to convert this series
\[
	1 - 1x + 2x^2 - 6x^3 + 24x^4 - 120x^5 + 720x^6 - 5040x^7 + \etc.
\]
into a continued fraction, extends further and can in the same way be applied to this much more general series
\begin{align*}
	z &= 1 -mx +m(m+n)x^2 - m(m+n)(m+2n)x^3 \\
	&+m(m+n)(m+2n)(m+3n)x^4 - \etc.;
\end{align*}
Then, having done the same operations, one will find
\[
	z = \cfrac{1}{1 + \cfrac{mx}{1 + \cfrac{nx}{1 + \cfrac{(m+n)x}{1 + \cfrac{2nx}{1 + \cfrac{(m+2n)x}{1 + \cfrac{3nx}{1 + \cfrac{(m+3n)x}{1 + \cfrac{4nx}{1 + \cfrac{(m+4n)x}{1 + \cfrac{5nx}{1 + \etc.}}}}}}}}}}}.
\]
But the same expression and other similar one can easily be found by means of the theorems, I proved in my dissertations on continued fractions in \textit{Comment. Acad. Petropol.}. Then I showed, that this equation
\[
	ax^{m-1}\mathrm{d}x = \mathrm{d}z + cx^{n-m-1}z\mathrm{d}x + bx^{n-1}z\mathrm{d}x
\]
is satisfied by this value
\[
	z = \cfrac{ax^m}{m + \cfrac{(ac+mb)x^n}{m + n + \cfrac{(ac-nb)x^n}{m + 2n + \cfrac{(ac + (m+n)b)x^n}{m + 3n + \cfrac{(ac-2nb)x^n}{m + 4n + \cfrac{(ac + (m+2n)b)x^n}{m + 5n + \cfrac{(ac-3nb)x^n}{m + 6n + \etc}.}}}}}}.
\]
Hence if $c=0$, it will be
\[
	\mathrm{d}z + bx^{n-1}z\mathrm{d}x = ax^{m-1}\mathrm{d}x
\]
and
\[
	e^{bx^n : n}z = a\int{e^{bx^n : n}x^{m-1}\mathrm{d}x} \quad \text{und} \quad z = ae^{-bx^n : n}\int{e^{bx^n : n}x^{m-1}\mathrm{d}x}
\]
and by a series
\[
	z = \frac{ax^m}{m} - \frac{abx^{m+n}}{m(m+n)} + \frac{ab^2 x^{m+2n}}{m(m+n)(m+2n)} - \frac{ab^3 x^{m+3n}}{m(m+n)(m+2n)(m+3n)} + \etc.
\]
But in this form our one we are treating is not contained.
\paragraph{§27}
But I further found, if one has this equation
\[
	fx^{m+n}\mathrm{d}x = x^{m+1}\mathrm{d}z + ax^{m}z\mathrm{d}x + bx^n z\mathrm{d}x + czz\mathrm{d}x,
\]
that the value of $z$ is expressed by a continued fraction of this kind
\[
	z = \cfrac{fx^m}{b + \cfrac{(mb+ab+cf)x^{m-n}}{b + \cfrac{(mb - nb + cf)x^{m-n}}{b + \cfrac{(2mb-nb+ab+cf)x^{m-n}}{b + \cfrac{(2mb-2nb+cf)x^{m-n}}{b + \cfrac{(3mb-2nb+ab+cf)x^{m-n}}{b + \cfrac{(3mb-3nb+cf)x^{m-n}}{b + \etc.}}}}}}}.
\]
Hence to be able to express the same value $z$ in a convenient way by  means of an ordinary series, let $c=0$, that one has this equation
\[
	fx^{m+n}\mathrm{d}x = x^{m+1}\mathrm{d}z + ax^m z\mathrm{d}x + bx^n z\mathrm{d}x,
\]
and by means of a continued fraction it will be
\[
	z = \cfrac{fx^m}{b + \cfrac{b(m+a)x^{m-n}}{b + \cfrac{b(m-n)x^{m-n}}{b + \cfrac{b(2m-n+a)x^{m-n}}{b + \cfrac{b(2m-2n)x^{m-n}}{b + \cfrac{b(3m-2n+a)x^{m-n}}{b + \cfrac{b(3m-3n)x^{m-n}}{b + \etc.}}}}}}}.
\]
By integration it will indeed be
\[
	x^a e^{bx^{n-m}:(n-m)}z = f\int{e^{bx^{n-m}:(m-n)}}x^{a+n-1}\mathrm{d}x
\]
or, if $m-n=k$, it will be
\[
	z = fe^{b:kx^k}x^{-a}\int{e^{-b:kx^k} x^{a+n-1}\mathrm{d}x},
\]
if one integrates in such a way of course, that $z$ vanishes for $x=0$. But by an infinite series it will be
\begin{align*}
	z =& \frac{f}{b}x^m - \frac{(m+a)}{b^2}fx^{2m-n} + \frac{(m+a)(2m-n+a)f}{b^3} x^{3m-2n} \\
	&-\frac{(m+a)(2m-n+a)(3m-2n+a)f}{b^4}x^{4m-3n} \\
	&+\frac{(m+a)(2m-n+a)(3m-2n+a)(4m-3n+a)f}{b^5}x^{5m-4n} - \etc.
\end{align*}
\paragraph{§28}
To simplify these expressions and at the same moment not restrict their generality, let us set
\[
	b = 1, \quad f = 1, \quad m+a=p, \quad m-n = q,
\]
that it is
\[
	a = p-m \quad \text{and} \quad n = m-q;
\]
and one will have this differential equation
\[
	x^m\mathrm{d}x = x^{q+1}\mathrm{d}z + (p-m)x^q z\mathrm{d}x + z\mathrm{d}x,
\]
whose integral is at first
\[
	z = e^{1:qx^q}x^{m-p}\int{e^{-18qx^q}x^{p-q-1}\mathrm{d}x}.
\]
The same value of the quantity $z$ will further be expressed by the following infinite series.
\[
	z = x^m - px^{m+q} + p(p+q)x^{m+2q} - p(p+q)(p+2q)x^{m+3q} + \etc.
\]
Finally this continued fraction will be equivalent to this series
\[
	z = \cfrac{x^m}{1 + \cfrac{px^q}{1 + \cfrac{qx^q}{1 + \cfrac{(p+q)x^q}{1 + \cfrac{2qx^q}{1 + \cfrac{(p+2q)x^q}{1 + \cfrac{3qx^q}{1 + \cfrac{(p+3q)x^q}{1 + \etc.}}}}}}}},
\]
which expression fully agrees with that, we obtained earlier in §$26$, and because there could be some doubt about the method, by which we found it, whether the numerators proceed according to the observed law to infinity or not, this doubt is now completely removed. Hence this consideration provided us with a method to sum innumerable divergent series or to find values equivalent to the same; among those that one, we treated, is a special case.
\paragraph{§29}
But further the case, in which $p=1$ and $q=2$ and $m=1$, seems be worthy to be noted; hence it will be
\[
	z = e^{1:2xx}\int{e^{-1:2xx}}\mathrm{d}x : xx
\]
and the infinite series will behave like this
\[
	z = x - 1x^3 + 1\cdot 3x^5 - 1\cdot 3\cdot 5x^7 + 1\cdot 3\cdot 5\cdot 7x^9 - \etc.,
\]
which is equal to this continued fraction
\[
	z = \cfrac{x}{1 + \cfrac{1xx}{1 + \cfrac{2xx}{1 + \cfrac{3xx}{1 + \cfrac{4xx}{1 + \cfrac{5xx}{1 + \cfrac{6xx}{1 + \etc.}}}}}}}.
\]
If therefore $x$ is set $=1$, that it is
\[
	z = 1 - 1 + 1\cdot 3 - 1\cdot 3\cdot 5 + 1\cdot 3\cdot 5\cdot 7 - 1\cdot 3\cdot 5\cdot 7\cdot 9 + \etc.,
\]
which series is strongly divergent, its value can nevertheless be expressed by this convergent continued fraction
\[
	z = \cfrac{1}{1 + \cfrac{1}{1 + \cfrac{2}{1 + \cfrac{3}{1 + \cfrac{4}{1 + \cfrac{5}{1 + \etc}}}}}}
\]
which yields the following fraction, approximately equal to the true value of $z$,

\begin{alignat*}{19}
	& &&  1 \quad && 2 \quad && 3 \quad && 4 \quad && ~5 \quad && ~6 \quad && ~7 \quad && ~~8 \quad && ~~9  \quad && ~~10  \quad && ~~11 \quad && ~~~12 & \\
	& z = && \frac{0}{1},~~ && \frac{1}{1},~~ && \frac{1}{2},~~ && \frac{3}{4},~~ && \frac{6}{10},~~ && \frac{18}{26},~~ && \frac{48}{76},~~ && \frac{156}{232},~~ &&\frac{492}{764},~~ && \frac{1740}{2620},~~ && \frac{6168}{9496},~~ && \frac{23568}{35696} ~~ && \etc.;
\end{alignat*}
hence if it is
\[
	z = \cfrac{1}{1 + \cfrac{1}{1 + \cfrac{2}{1 + \cfrac{3}{1 + \cfrac{4}{1 + \cfrac{5}{1 + \cfrac{6}{1 + \cfrac{7}{1 + \cfrac{8}{1 + \cfrac{9}{1 + \cfrac{10}{1 + p}}}}}}}}}}},
\]
it will be
\[
	z = \frac{23568 + 6168p}{35696 + 9496p}
\]
or
\[
	z = \frac{2946 + 771p}{4402 + 1187p}
\]
and
\[
	p = \cfrac{11}{1 + \cfrac{12}{1 + \cfrac{13}{1 + \cfrac{14}{1 + \cfrac{15}{1 + \etc.}}}}}
\]
Let be
\[
	p = \frac{11}{1+q} \quad \text{and} \quad q = \frac{12}{1 + r}
\]
it will be
\[
r=\frac{12-q}{q}
\]
and since $p$, $q$, $r$ grow uniformly, it will be
\[
	2q = \frac{12 + 22q -qq}{q+qq} \quad \text{and} \quad 2q^3 + 3qq - 22q - 12 = 0,
\]
where it is approximately
\[
	q = 2,94,\quad p = 2,79 \quad \text{and} \quad z = \frac{5097,09}{7773,73} = 0,65568.
\]
\end{document}